\documentclass[a4paper,12pt]{article}
 
 \usepackage{latexsym,amssymb,amsmath}  
 
 \newtheorem{thm}{Theorem}
 \newtheorem{lem}{Lemma}

 \newtheorem{defn}{Definition} 

\begin{document}

\title{Fueter's theorem for the biregular functions of Clifford analysis}

\author{Dixan Pe\~na Pe\~na$^{*}$ and Frank Sommen$^{**}$}

\date{\normalsize{Clifford Research Group, Department of Mathematical Analysis\\Faculty of Engineering, Ghent University\\Galglaan 2, 9000 Gent, Belgium}\\\vspace{0.1cm}
\small{$^{*}$e-mail: dpp@cage.ugent.be\\
$^{**}$e-mail: fs@cage.ugent.be}}

\maketitle

\begin{abstract}
\noindent In this paper we present a generalization of the Fueter's theorem for monogenic functions to the case of the biregular functions.\vspace{0.2cm}\\
\textit{Keywords}: Clifford algebras; biregular functions; Fueter's theorem.\vspace{0.1cm}\\
\textit{Mathematics Subject Classification}: 30G35.
\end{abstract}

\section{Introduction}

Let $\mathbb{R}_{0,m}$ be the $2^m$-dimensional real Clifford algebra generated by the standard basis $\{e_1,\ldots,e_m\}$ of the Euclidean space $\mathbb R^m$ (see \cite{Cl}). The multiplication in $\mathbb{R}_{0,m}$ is determined by the relations 
\[e_je_k+e_ke_j=-2\delta_{jk},\quad j,k=1,\dots,m\]
where $\delta_{jk}$ denotes the Kronecker delta. A basis for the algebra is then given by the elements $e_A=e_{j_1}\dots e_{j_k}$ where $A=\{j_1,\dots,j_k\}\subset\{1,\dots,m\}$ and $j_1<\dots<j_k$ ($e_{\emptyset}=1$ is the identity element). A general element $a$ of $\mathbb R_{0,m}$ may thus be written as
\[a=\sum_Aa_Ae_A,\quad a_A\in\mathbb R\]
and its conjugate $\overline a$ is defined by
\[\overline a=\sum_Aa_A\overline e_A,\quad\overline e_A=\overline e_{j_k}\dots\overline e_{j_1},\quad\overline e_j=-e_j,\;j=1,\dots,m.\]
Suppose now that $\Omega$ is an open subset of $\mathbb R^{m+1}\times\mathbb R^{m+1}$ and let $f$ be an $\mathbb{R}_{0,m}$-valued function defined in $\Omega$. Then $f$ is of the form 
\[f(x,y)=\sum_{A}f_A(x,y)e_A,\quad(x,y)=(x_0,\dots,x_m,y_0,\dots,y_m)\in\Omega\]
where the functions $f_A$ are $\mathbb R$-valued. The variables $x$ and $y$ will also be identified with the paravectors $x=x_0+\underline x=x_0+\sum_{j=1}^mx_je_j$ and $y=y_0+\underline y=y_0+\sum_{j=1}^my_je_j$, respectively. Next, we introduce the generalized Cauchy-Riemann operators
\[\partial_{x}=\partial_{x_0}+\partial_{\underline x}=\partial_{x_0}+\sum_{j=1}^me_j\partial_{x_j},\quad\partial_y=\partial_{y_0}+\partial_{\underline y}=\partial_{y_0}+\sum_{j=1}^me_j\partial_{y_j},\]
which factorize the Laplacian in the variables $x$ and $y$, respectively, i.e.
\[\Delta_x=\sum_{j=0}^m\partial_{x_j}^2=\partial_x\overline\partial_x=\overline\partial_x\partial_x,\quad\Delta_y=\sum_{j=0}^m\partial_{y_j}^2=\partial_y\overline\partial_y=\overline\partial_y\partial_y.\]
In this paper we shall deal with the so-called biregular functions, which are defined as follows. 

\begin{defn}
Assume that $f\in C^1(\Omega)$. Then the function $f$ is called biregular in $\Omega$ if it fulfills in $\Omega$ the system $\partial_xf=f\partial_y=0$.
\end{defn}

The biregular functions were introduced in the 1980s by Brackx and Pincket as an extension to two higher dimensional variables of the standard monogenic functions, i.e. $C^1$ functions $f:\Omega\subset\mathbb R^{m+1}\rightarrow\mathbb{R}_{0,m}$ satisfying $\partial_xf=0$ (or $f\partial_x=0$). Some of the main properties of the biregular functions may be found in e.g. \cite{BP1,BP2,BP3,S1,S2}. For a detailed study of the monogenic functions we refer the reader to \cite{BDS,DSS,GuSp}.

An important technique to generate monogenic functions is the so-called Fueter's theorem. Discovered by R. Fueter in the setting of quaternionic analysis (see \cite{F}), this technique has been extended to $\mathbb R^{m+1}$ within the framework of Clifford analysis in \cite{Sce,S} ($m$ odd) and in \cite{Q} ($m$ even). For other works on this topic we refer the reader to e.g. \cite{CoSaF,CoSaF2,KQS,DS1,DS2,DS3,QS}.

For $m$ odd, Fueter's theorem states the following: 

\begin{thm}\label{FueThmF}
Let $u+iv$ be a holomorphic function in some open subset $\Xi$ of the upper half of the complex plane $\mathbb C$ and assume that $P_k(\underline x)$ is a homogeneous monogenic polynomial of degree $k$ in $\mathbb R^m$. Put $\underline\omega=\underline x/r$, with $r=\vert\underline x\vert$. If $m$ is odd, then the function
\begin{equation*}
\Delta_x^{k+\frac{m-1}{2}}\bigl[\bigl(u(x_0,r)+\underline\omega\,v(x_0,r)\bigr)P_k(\underline x)\bigr]
\end{equation*}
is (left) monogenic in $\Omega=\{x\in\mathbb R^{m+1}:\;(x_0,r)\in\Xi\}$. 
\end{thm}

Fueter's theorem discloses a remarkable connection existing between the classical holomorphic functions and its higher dimensional counterpart (i.e. the monogenic functions).  The purpose of this paper is to generalize this important result to the case of the biregular functions. We shall see how a similar relationship exists between the latter and the holomorphic functions of two complex variables.

\section{Auxiliary results}

We group together here a series of lemmas which will be needed in the proof of our main result. 

For an $\mathbb R$-valued $C^1$ function $\phi$ and an $\mathbb R_{0,m}$-valued $C^1$ function $g$, it is clear that 
\begin{equation}\label{Leibnitz1}
\partial_{\underline x}(\phi g)=(\partial_{\underline x}\phi)g+\phi(\partial_{\underline x}g),
\end{equation}
\begin{equation}\label{Leibnitz2}
(\phi g)\partial_{\underline y}=g(\partial_{\underline y}\phi)+\phi(g\partial_{\underline y}).
\end{equation}
Moreover, for a vector-valued $C^1$ function $\underline f=\sum_{j=1}^mf_je_j$, we have 
\begin{equation}\label{lr1}
\partial_{\underline x}(\underline f\,g)=(\partial_{\underline x}\underline f)\,g-\underline f(\partial_{\underline x}g)-2\sum_{j=1}^mf_j(\partial_{x_j}g).
\end{equation}
Indeed,
\begin{equation*}
\partial_{\underline x}(\underline fg)=\sum_{j=1}^me_j\left((\partial_{x_j}\underline f)g+\underline f(\partial_{x_j}g)\right)=(\partial_{\underline x}\underline f)g+\sum_{j=1}^me_j\underline f(\partial_{x_j}g),
\end{equation*}
which results in (\ref{lr1}) on account of the equality 
\[e_j\underline f=-\underline fe_j-2f_j,\quad j=1,\dots,m.\]
In the same spirit we can also prove:
\begin{equation}\label{lr2}
(g\underline f)\partial_{\underline y}=g(\underline f\partial_{\underline y})-(g\partial_{\underline y})\underline f-2\sum_{j=1}^mf_j(\partial_{y_j}g).
\end{equation}

\begin{lem}\label{lem1}
Suppose that $f(t_1,\dots,t_d)$ is an $\mathbb R$-valued $C^\infty$ function on $\mathbb R^d$ and that $D_{t_j}$ and $D^{t_j}$ are differential operators defined by 
\begin{align*}
D_{t_j}(n)\{f\}&=\left(\frac{1}{{t_j}}\,\partial_{t_j}\right)^nf,\quad j=1,\dots,d,\\
D^{t_j}(n)\{f\}&=\partial_{t_j}\left(\frac{D^{t_j}(n-1)\{f\}}{{t_j}}\right),\quad j=1,\dots,d,
\end{align*}
for $n\ge1$ and $D_{t_j}(0)\{f\}=D^{t_j}(0)\{f\}=f$. Then one has
\begin{itemize}
\item[{\rm(i)}] $\partial_{t_j}^2D_{t_j}(n)\{f\}=D_{t_j}(n)\{\partial_{t_j}^2f\}-2n D_{t_j}(n+1)\{f\}$, 
\item[{\rm(ii)}] $\partial_{t_j}D_{t_j}(n-1)\left\{f/{t_j}\right\}=D^{t_j}(n)\{f\}$,
\item[{\rm(iii)}] $D^{t_j}(n)\{\partial_{t_j}f\}=\partial_{t_j}D_{t_j}(n)\{f\}$,  
\item[{\rm(iv)}] $D_{t_j}(n)\{\partial_{t_j}f\}-\partial_{t_j}D^{t_j}(n)\{f\}=2n/t_j\,D^{t_j}(n)\{f\}$,
\item[{\rm(v)}] $\partial_{t_j}^2D^{t_j}(n)\{f\}=D^{t_j}(n)\{\partial_{t_j}^2f\}-2nD^{t_j}(n+1)\{f\}$.
\end{itemize}
\end{lem}
\textit{Proof.} We prove (i) by induction. When $n=1$, we have  
\begin{align*}
\partial_{t_j}^2D_{t_j}(1)\{f\}&=\frac{\partial_{t_j}^3f}{t_j}-2\frac{\partial_{t_j}^2f}{t_j^2}+2\frac{\partial_{t_j}f}{t_j^3}\\&=D_{t_j}(1)\big\{\partial_{t_j}^2f\big\}-2D_{t_j}(2)\{f\}
\end{align*}
as desired.

Now we proceed to show that when (i) holds for a positive integer $n$, then it also holds for $n+1$. Indeed,
\begin{align*}
\partial_{t_j}^2D_{t_j}(n+1)\{f\}&=D_{t_j}(1)\big\{\partial_{t_j}^2D_{t_j}(n)\{f\}\big\}-2D_{t_j}(2)\big\{D_{t_j}(n)\{f\}\big\}\\&=D_{t_j}(1)\Bigl\{D_{t_j}(n)\big\{\partial_{t_j}^2f\big\}-2n\,D_{t_j}(n+1)\{f\}\Bigr\}\\&\quad-2D_{t_j}(n+2)\{f\}\\&=D_{t_j}(n+1)\big\{\partial_{t_j}^2f\big\}-2(n+1)\,D_{t_j}(n+2)\{f\}.
\end{align*}
Statement (ii) easily follows from the definition of $D^{t_j}(n)\{f\}$. Next, using (ii), we obtain (iii) as
\[D^{t_j}(n)\{\partial_{t_j}f\}=\partial_{t_j}D_{t_j}(n-1)\{\partial_{t_j}f/t_j\}=\partial_{t_j}D_{t_j}(n)\{f\}.\]
To obtain (iv) we use (i) and (ii):
\begin{align*}
D_{t_j}(n)&\{\partial_{t_j}f\}-\partial_{t_j}D^{t_j}(n)\{f\}\\
&=D_{t_j}(n)\{\partial_{t_j}f\}-\partial_{t_j}^2D_{t_j}(n-1)\{f/t_j\}\\
&=D_{t_j}(n)\{\partial_{t_j}f\}-D_{t_j}(n-1)\big\{\partial_{t_j}^2\{f/t_j\}\big\}+2(n-1)\,D_{t_j}(n)\{f/t_j\}\\
&=D_{t_j}(n)\{\partial_{t_j}f\}-D_{t_j}(n-1)\big\{D_{t_j}(1)\{\partial_{t_j}f\}-2D_{t_j}(1)\{f/t_j\}\big\}\\
&\qquad+2(n-1)\,D_{t_j}(n)\{f/t_j\}\\
&=2n\,D_{t_j}(n)\{f/t_j\}=2n/t_j\,D^{t_j}(n)\{f\}.
\end{align*}
Finally, from (i)-(iii) it follows that
\begin{align*}
\partial_{t_j}^2D^{t_j}(n)\{f\}&=\partial_{t_j}^3D_{t_j}(n-1)\{f/t_j\}\\
&=\partial_{t_j}D_{t_j}(n-1)\{\partial_{t_j}^2\{f/t_j\}\}-2(n-1)\partial_{t_j}D_{t_j}(n)\{f/t_j\}\\
&=\partial_{t_j}D_{t_j}(n)\{\partial_{t_j}f\}-2n\partial_{t_j}D_{t_j}(n)\{f/t_j\}\\
&=D^{t_j}(n)\{\partial_{t_j}^2f\}-2nD^{t_j}(n+1)\{f\},
\end{align*}
thus proving (v).\hfill$\square$\vspace{0.2cm}

Throughout this paper we denote by $P_{k,l}(\underline x,\underline y)$ a homogeneous biregular polynomial in $\mathbb R^m\times\mathbb R^m$ of degree $k$ in $\underline x$ and degree $l$ in $\underline y$, i.e.
\begin{alignat*}{2}
\partial_{\underline x}P_{k,l}(\underline x,\underline y)&=P_{k,l}(\underline x,\underline y)\partial_{\underline y}=0,&\qquad(\underline x,\underline y)&\in\mathbb R^m\times\mathbb R^m,\\
P_{k,l}(t_1\underline x,t_2\underline y)&=t_1^kt_2^lP_{k,l}(\underline x,\underline y),&\qquad t_1,t_2&\in\mathbb R.
\end{alignat*}
Put 
\[\underline\omega=\underline x/r,\quad\underline\nu=\underline y/\rho\]
where $r=\vert\underline x\vert$ and $\rho=\vert\underline y\vert$. 

\begin{lem}\label{lem2}
Let $h(x_0,r,y_0,\rho)$ be an $\mathbb R$-valued $C^\infty$ function on $\mathbb R^4$ such that \[\partial_{x_0}^2h+\partial_{r}^2h=\partial_{y_0}^2h+\partial_{\rho}^2h=0.\]
Then
\[\Delta_x^n\big(hP_{k,l}\big)=\prod_{j=1}^n(2k+m-(2j-1))D_r(n)\{h\}P_{k,l},\]
\[\Delta_x^n(h\,\underline\omega\,P_{k,l})=\prod_{j=1}^n(2k+m-(2j-1))D^r(n)\{h\}\underline\omega P_{k,l},\]
\[\Delta_y^n(hP_{k,l})=\prod_{j=1}^n(2l+m-(2j-1))D_\rho(n)\{h\}P_{k,l},\]
\[\Delta_y^n(hP_{k,l}\,\underline\nu)=\prod_{j=1}^n(2l+m-(2j-1))D^\rho(n)\{h\}P_{k,l}\,\underline\nu.\]
\end{lem}
\textit{Proof.} We first prove that for any $\mathbb R$-valued $C^2$ function $g(x_0,r,y_0,\rho)$ in the variables $x_0$, $r$, $y_0$ and $\rho$ the following equalities hold  
\[\Delta_x(gP_{k,l})=\left(\partial_{x_0}^2g+\partial_r^2g+(2k+m-1)D_r(1)\{g\}\right)P_{k,l},\]
\[\Delta_x(g\,\underline\omega\,P_{k,l})=\left(\partial_{x_0}^2g+\partial_r^2g+(2k+m-1)D^r(1)\{g\}\right)\underline\omega\,P_{k,l},\]
\[\Delta_y(gP_{k,l})=\left(\partial_{y_0}^2g+\partial_\rho^2g+(2l+m-1)D_\rho(1)\{g\}\right)P_{k,l},\]
\[\Delta_y(gP_{k,l}\,\underline\nu)=\left(\partial_{y_0}^2g+\partial_\rho^2g+(2l+m-1)D^\rho(1)\{g\}\right)P_{k,l}\,\underline\nu.\]
In fact, it follows that
\[\partial_{\underline x}g=\sum_{j=1}^me_j\partial_{x_j}g=\sum_{j=1}^me_j(\partial_rg)(\partial_{x_j}r)=\underline\omega\partial_rg,\]
and
\[\Delta_{\underline x}\,\underline\omega=-\partial_{\underline x}^2\,\underline\omega=(m-1)\partial_{\underline x}\left(\frac{1}{r}\right)=-\frac{(m-1)}{r^2}\,\underline\omega,\]
\begin{align*}
\Delta_xg&=\partial_{x_0}^2g+\Delta_{\underline x}g=\partial_{x_0}^2g-\partial_{\underline x}(\underline\omega\partial_rg)\\
&=\partial_{x_0}^2g+\partial_r^2g+\frac{m-1}{r}\,\partial_rg.
\end{align*}
Therefore 
\begin{align*}
\Delta_x(gP_{k,l})&=(\Delta_xg)P_k+g(\Delta_{\underline x}P_{k,l})+2\sum_{j=1}^m(\partial_{x_j}g)(\partial_{x_j}P_{k,l})\\
&=\left(\partial_{x_0}^2g+\partial_r^2g+\frac{m-1}{r}\,\partial_rg\right)P_{k,l}+2\frac{\partial_rg}{r}\sum_{j=1}^mx_j(\partial_{x_j}P_{k,l})\\
&=\left(\partial_{x_0}^2g+\partial_r^2g+\frac{2k+m-1}{r}\,\partial_rg\right)P_{k,l}
\end{align*}
where we have also used Euler's theorem for homogeneous functions. Moreover,
\begin{align*}
\Delta_x(g\underline\omega P_{k,l})&=(\Delta_{\underline x}\,\underline\omega)gP_{k,l}+\underline\omega\Delta_x(gP_{k,l})+2\sum_{j=1}^m(\partial_{x_j}\underline\omega)(\partial_{x_j}(gP_{k,l}))\\
&=-\frac{(m-1)}{r^2}\,g\underline\omega P_{k,l}+\underline\omega\Delta_x(gP_{k,l})\\
&\qquad+2\sum_{j=1}^m\left(\frac{e_j}{r}-\frac{x_j}{r^2}\,\underline\omega\right)\left(\frac{x_j}{r}\,(\partial_rg)P_{k,l}+g(\partial_{x_j}P_{k,l})\right)\\
&=\left(\partial_{x_0}^2g+\partial_r^2g+(2k+m-1)\left(\frac{\partial_rg}{r}-\frac{g}{r^2}\right)\right)\underline\omega P_{k,l}
\end{align*}
In a similar way we can prove the other two equalities. 

The proof now follows by induction using the previous equalities together with statements (i) and (v) of Lemma \ref{lem1}. 

It is clear that the lemma is true in the case $n=1$. Assume that the formulae hold for a positive integer $n$; we will prove them for $n+1$. 

We thus get
\[\Delta_x^{n+1}\big(hP_{k,l}\big)=\prod_{j=1}^n\big(2k+m-(2j-1)\big)\Delta_x\big(D_r(n)\{h\}P_{k,l}\big).\]
But
\begin{align*}
&\Delta_x\big(D_r(n)\{h\}P_{k,l}\big)\\
&=\left(\partial_{x_0}^2D_r(n)\{h\}+\partial_r^2D_r(n)\{h\}+(2k+m-1)D_r(n+1)\{h\}\right)P_{k,l}\\
&=\left(D_r(n)\{\partial_{x_0}^2h+\partial_r^2h\}+\big(2k+m-(2n+1)\big)D_r(n+1)\{h\}\right)P_{k,l}
\end{align*}
yielding
\[\Delta_x^{n+1}\big(hP_{k,l}\big)=\prod_{j=1}^{n+1}\big(2k+m-(2j-1)\big)D_r(n+1)\{h\}P_{k,l},\]
which establishes the first formula. The others may be proved similarly.\hfill$\square$

\begin{lem}\label{lem3}
Assume that $A$, $B$, $C$ and $D$ are $\mathbb R$-valued $C^1$ functions on $\mathbb R^4$. Then the function  
\begin{multline*}
F(x,y)=A(x_0,r,y_0,\rho)P_{k,l}(\underline x,\underline y)+B(x_0,r,y_0,\rho)\,\underline\omega\,P_{k,l}(\underline x,\underline y)\\
+C(x_0,r,y_0,\rho)P_{k,l}(\underline x,\underline y)\,\underline\nu+D(x_0,r,y_0,\rho)\,\underline\omega\,P_{k,l}(\underline x,\underline y)\,\underline\nu,
\end{multline*}
is biregular if the following Vekua-type systems are satisfied
\begin{align*}
\partial_{x_0}A-\partial_rB&=\frac{2k+m-1}{r}\,B & \partial_{x_0}C-\partial_rD&=\frac{2k+m-1}{r}\,D\\
\partial_{x_0}B+\partial_rA&=0 & \partial_{x_0}D+\partial_rC&=0
\end{align*}
\begin{align*}
\partial_{y_0}A-\partial_\rho C&=\frac{2l+m-1}{\rho}\,C & \partial_{y_0}B-\partial_\rho D&=\frac{2l+m-1}{\rho}\,D\\ 
\partial_{y_0}C+\partial_\rho A&=0 & \partial_{y_0}D+\partial_\rho B&=0.
\end{align*}
\end{lem}
\textit{Proof.} It is not hard to check, using the Leibniz rules $(\ref{Leibnitz1})$, $(\ref{lr1})$ and Euler's theorem for homogeneous functions, that the following equality holds
\begin{multline*}
\partial_{\underline x}F(x,y)=\partial_rA\,\underline\omega\,P_{k,l}-\left(\partial_rB+\frac{2k+m-1}{r}\,B\right)P_{k,l}+\partial_rC\,\underline\omega\,P_{k,l}\,\underline\nu\\
-\left(\partial_rD+\frac{2k+m-1}{r}\,D\right)P_{k,l}\,\underline\nu.
\end{multline*}
Therefore
\begin{multline*}
\partial_xF(x,y)=\left(\partial_{x_0}A-\partial_rB-\frac{2k+m-1}{r}\,B\right)P_{k,l}+\left(\partial_{x_0}B+\partial_rA\right)\,\underline\omega\,P_{k,l}\\
\left(\partial_{x_0}C-\partial_rD-\frac{2k+m-1}{r}\,D\right)P_{k,l}\,\underline\nu+(\partial_{x_0}D+\partial_rC)\,\underline\omega\,P_{k,l}\,\underline\nu.
\end{multline*}
Similarly, using $(\ref{Leibnitz2})$, $(\ref{lr2})$ and Euler's theorem for homogeneous functions, we can also obtain that
\begin{align*}  
\begin{split}
F(x,y)&\partial_y=\left(\partial_{y_0}A-\partial_\rho C-\frac{2l+m-1}{\rho}\,C\right)P_{k,l}\\
&\qquad+\left(\partial_{y_0}B-\partial_\rho D-\frac{2l+m-1}{\rho}\,D\right)\,\underline\omega\,P_{k,l}+(\partial_{y_0}C+\partial_\rho A)\,P_{k,l}\,\underline\nu\\
&\qquad+(\partial_{y_0}D+\partial_\rho B)\,\underline\omega\,P_{k,l}\,\underline\nu,
\end{split}
\end{align*}
which completes the proof.\hfill$\square$

\section{Fueter's theorem}

Suppose that $u_j(x_1,y_1,x_2,y_2)$ and $v_j(x_1,y_1,x_2,y_2)$, $j=1,2$, are $\mathbb R$-valued $C^1$ functions which satisfy the following Cauchy-Riemann systems 
\begin{align*}
\partial_{x_1}u_1&=\partial_{y_1}v_1 & \partial_{x_1}u_2&=\partial_{y_1}v_2\\ 
\partial_{y_1}u_1&=-\partial_{x_1}v_1 & \partial_{y_1}u_2&=-\partial_{x_1}v_2
\end{align*}
\begin{align*}
\partial_{x_2}u_1&=\partial_{y_2}u_2&\partial_{x_2}v_1&=\partial_{y_2}v_2\\
\partial_{y_2}u_1&=-\partial_{x_2}u_2&\partial_{y_2}v_1&=-\partial_{x_2}v_2
\end{align*} 
in some open subset $\Xi\subset\mathbb R^2_+\times\mathbb R^2_+=\{(x_1,y_1,x_2,y_2)\in\mathbb R^2\times\mathbb R^2:\;y_1,\,y_2>0\}$. In other words, the functions $u_1+iv_1$ and $u_2+iv_2$ are holomorphic with respect to the complex variable $z_1=x_1+iy_1$ while $u_1+iu_2$ and $v_1+iv_2$ are holomorphic with respect to $z_2=x_2+iy_2$.

\noindent\textbf{Remark:} Note that a simple way to obtain these types of functions is by using holomorphic functions of two complex variables. Indeed, if $u+iv$ is a holomorphic function in $\mathbb C^2$, then we can put $u_1=u$, $v_1=u_2=v$ and $v_2=-u$.    

Finally, for the biregular functions we propose the following result.

\begin{thm}
Let $u_j$, $v_j$ $(j=1,2)$ be as above. If $m$ is odd, then the function
\begin{multline*}
\Delta_x^{k+\frac{m-1}{2}}\Delta_y^{l+\frac{m-1}{2}}\biggl(u_1(x_0,r,y_0,\rho)P_{k,l}(\underline x,\underline y)+v_1(x_0,r,y_0,\rho)\,\underline\omega\,P_{k,l}(\underline x,\underline y)\\
+u_2(x_0,r,y_0,\rho)P_{k,l}(\underline x,\underline y)\,\underline\nu+v_2(x_0,r,y_0,\rho)\,\underline\omega\,P_{k,l}(\underline x,\underline y)\,\underline\nu\biggr)
\end{multline*}
is biregular in $\Omega=\{(x,y)\in\mathbb R^{m+1}\times\mathbb R^{m+1}:\;(x_0,r,y_0,\rho)\in\Xi\}$.
\end{thm}
\textit{Proof.} By Lemma \ref{lem2}, we get that 
\begin{multline*}
\Delta_x^{k+\frac{m-1}{2}}\Delta_y^{l+\frac{m-1}{2}}\Big(u_1P_{k,l}+v_1\,\underline\omega\,P_{k,l}+u_2P_{k,l}\,\underline\nu+v_2\,\underline\omega\,P_{k,l}\,\underline\nu\Big)\\
=(2k+m-1)!!(2l+m-1)!!\Big(AP_{k,l}+B\,\underline\omega\,P_{k,l}+CP_{k,l}\,\underline\nu+D\,\underline\omega\,P_{k,l}\,\underline\nu\Big),
\end{multline*}
with 
\[A=D_r\left(k+\frac{m-1}{2}\right)D_\rho\left(l+\frac{m-1}{2}\right)\{u_1\},\]
\[B=D^r\left(k+\frac{m-1}{2}\right)D_\rho\left(l+\frac{m-1}{2}\right)\{v_1\},\]
\[C=D_r\left(k+\frac{m-1}{2}\right)D^\rho\left(l+\frac{m-1}{2}\right)\{u_2\},\]
\[D=D^r\left(k+\frac{m-1}{2}\right)D^\rho\left(l+\frac{m-1}{2}\right)\{v_2\}.\]
The task is now to prove that $A$, $B$, $C$ and $D$ satisfy the Vekua-type systems of Lemma \ref{lem3}. In order to do that, it will be necessary to use the assumptions on $u_j$ and $v_j$ $(j=1,2)$ and Lemma \ref{lem1}. 

Indeed, by the statements (iii) and (iv) of Lemma \ref{lem1} and using the fact that $u_1+iv_1$ is holomorphic with respect to the complex variable $z_1=x_1+iy_1$, it follows that 
\[\begin{split}
&\partial_{x_0}A-\partial_rB\\
&= D_\rho\left(l+\frac{m-1}{2}\right)\left\{D_r\left(k+\frac{m-1}{2}\right)\{\partial_{x_0}u_1\}-\partial_rD^r\left(k+\frac{m-1}{2}\right)\{v_1\}\right\}\\
&=D_\rho\left(l+\frac{m-1}{2}\right)\left\{D_r\left(k+\frac{m-1}{2}\right)\{\partial_rv_1\}-\partial_rD^r\left(k+\frac{m-1}{2}\right)\{v_1\}\right\}\\
&=\frac{2k+m-1}{r}D^r\left(k+\frac{m-1}{2}\right)D_\rho\left(l+\frac{m-1}{2}\right)\{v_1\}\\
&=\frac{2k+m-1}{r}\,B
\end{split}\]
and 
\[\begin{split}
&\partial_{x_0}B+\partial_rA\\
&=D_\rho\left(l+\frac{m-1}{2}\right)\left\{D^r\left(k+\frac{m-1}{2}\right)\{\partial_{x_0}v_1\}+\partial_rD_r\left(k+\frac{m-1}{2}\right)\{u_1\}\right\}\\
&=D_\rho\left(l+\frac{m-1}{2}\right)\left\{D^r\left(k+\frac{m-1}{2}\right)\{\partial_{x_0}v_1\}+D^r\left(k+\frac{m-1}{2}\right)\{\partial_ru_1\}\right\}\\
&=D^r\left(k+\frac{m-1}{2}\right)D_\rho\left(l+\frac{m-1}{2}\right)\{\partial_{x_0}v_1+\partial_ru_1\}\\
&=0.
\end{split}\]
In a similar way we can verify the other systems of Lemma \ref{lem3}.\hfill$\square$

\subsection*{Acknowledgments}

D. Pe\~na Pe\~na acknowledges the support of a Postdoctoral Fellowship from \lq\lq Special Research Fund" (BOF) of Ghent University.

\end{document}